\documentclass{article}
\usepackage{amssymb}
\usepackage{paralist}
\usepackage{changes}
\usepackage{graphics} 
\usepackage{epsfig} 
\usepackage{graphicx}  \usepackage{epstopdf}
\usepackage[colorlinks=true]{hyperref}
\hypersetup{urlcolor=blue, citecolor=red}

\newtheorem{remark}{Remark}

\title{New insights on numerical error in symplectic integration
\thanks{The first author is supported by a grant from the Fondation du Coll\`ege de France under the research convention PU14150472.}
}

\newcommand{\bz}{{\bf z}}
\newcommand{\bZ}{{\bf Z}}
\newcommand{\bp}{{\bf p}}
\newcommand{\bq}{{\bf q}}
\newcommand{\bP}{{\bf P}}
\newcommand{\bQ}{{\bf Q}}
\newcommand{\bu}{{\bf u}}
\newcommand{\bv}{{\bf v}}

\begin{document}
\maketitle

\centerline{\scshape Hugo Jim\'enez-P\'erez$^*$, Jean-Pierre Vilotte}
\medskip
{\footnotesize
 \centerline{Institut de Physique du Globe de Paris}
   \centerline{1 rue Jussieu, 75865 Cedex Paris, France}
} 

\centerline{\scshape Barbara Romanowicz$^{1,2,3}$}
\medskip
{\footnotesize
 \centerline{$^1$ Institut de Physique du Globe de Paris}
   \centerline{1 rue Jussieu, 75865 Cedex Paris, France}
 \centerline{$^2$ Coll\`ege de France, Paris, France}
 \centerline{$^3$ Department of Earth and Planetary Sciences}
   \centerline{University of California, Berkeley, USA}
} 

\begin{abstract}
We implement and investigate the numerical properties of a new family of integrators 
connecting both variants of the symplectic Euler schemes, and including
an alternative to the classical symplectic mid-point scheme, with some additional terms. 
This family is derived from a new algorithm, introduced in a previous study,
for generating symplectic integrators based on the concept of special symplectic 
manifold. The use of symplectic rotations and a particular type of projection
keeps the whole procedure within the symplectic framework. 

We show that it is  possible to define a set of parameters that control the 
additional terms providing a way of ``tuning''
these new symplectic schemes. We test the  ``tuned'' symplectic integrators with the
perturbed pendulum and we compare its behavior with an explicit 
$\mathcal {SABA}_2$ scheme for perturbed systems. Remarkably, for the given examples,
the error in the energy integral can be reduced considerably. 
There is a natural geometrical explanation, sketched at the end of this paper.
This is the subject of a parallel article where a finer analysis is performed. 
Numerical results obtained in this paper open a new point of view on symplectic 
integrators and Hamiltonian error.
\end{abstract}

\section{Introduction}

A symplectic integrator for a Hamiltonian system is a numerical method 
which preserves the structure of the Hamiltonian vector field. Poincar\'e
discovered that the flow of a Hamiltonian system forms a one-parameter subgroup 
of canonical transformations. In modern language, we say that the Hamiltonian 
flow is symplectic. The standard procedure for simulating Hamiltonian dynamics is
by discretizing the Hamiltonian flow, which consists in discretizing the evolution 
time and looking for symplectic transformations which map the state of the system 
between two adjacent elements of the discretized time, e.g. from time $t_n$ to
time $t_{n+1}$.

It is well-known that one method for creating symplectic maps is based on generating 
functions, and it was already used by Poincar\'e when looking for periodic orbits 
of second genus \cite{Poi99}. In fact, generating functions were an important 
ingredient of the invariant integral theory pioneered by Poincar\'e 
and generalized by Cartan \cite{Car22}.
Since then, interest in generating functions remains very active from both the theoretical 
and the numerical point of view. Indeed for every symplectic 
transformation $\psi$ there corresponds (at least locally) a class of functions generating
a Lagrangian submanifold, i.e., the graph of a 1-form symplectomorphic to
the Liouville form. In addition, Hamilton-Jacobi theory connects this Lagrangian 
submanifold with another submanifold invariant under the transformation $\psi$.

Lagrangian submanifolds used to obtain suitable maps for
symplectic integrators must have a very particular form. They must 
contain all the information concerning the source and the target 
(symplectic) variables. This information, encoded in the Lagrangian submanifold,
contains the well-known fact that generating functions for 
symplectic integrators must generate maps close to the identity. This 
criterion is not enough for determining whether the generating function 
associated to some symplectic map is suitable for obtaining a symplectic integrator.
What is generally missing in the literature, is a geometrical approach 
distilling the theory behind the different techniques within a unified
point of view on symplectic integrators. The present paper, together with 
\cite{Jim15a, Jim15c} are contributions in this direction. Before going over 
the theory, we briefly review previous results
on symplectic integrators and generating functions related to the present paper.

The first article dealing with symplectic algorithms is attributed to 
De Vogelaere \cite{DeV56} in 1956. In 1983 Ruth \cite{Rut83} and Channell \cite{Cha83}
made some progress with different techniques, in particular, Ruth pioneered 
explicit symplectic integrators and composition-splitting methods. 
Additional contributions were made independently by Menyuk \cite{Men84} and 
Kang \cite{Kang85a} in 1984.
The same year, Kang Feng and his collaborators
started a systematic study of symplectic integrators 
using generating functions \cite{Kang85b,GM88,KHMD89}.
His point of view was mostly algebraic, and based on Siegel's article \cite{Sie43},
reprinted some years later in book format \cite{Sie64}. Some geometrization was 
achieved by  Ge and Marsden \cite{GM88}, Ge \cite{Ge90,Ge91}, 
Ge and Dau-liu \cite{GD95}, although their numerical algorithms were based on 
Feng's procedure\footnote{Most of Feng's articles and some from his collaborators 
were recently edited in book format by M. Qin \cite{KQ10}.}. 
In 1990 Channell and Scovel \cite{CS90}
introduced symbolic computations to derive the integration formulas which arose in the procedure.
Other important contributions were made by Sanz-Serna \cite{Sanz88} who worked on
symplecticity conditions for Runge-Kutta methods, and Miesbach and Pesch \cite{MP92}
who introduced some methods using Runge-Kutta techniques with generating functions.
Many other authors have produced algorithms using generating functions, but
the geometrical construction remains the same, and their contributions diverge from
our discussion.

In a recent work \cite{Jim15a}, another link between generating functions for
symplectic integrators and symplectic geometry has been studied,
based on the concept of \emph{special symplectic manifold} introduced by
Tulczjyew in \cite{Tul76, Tul77}. Generating forms and functions in this 
framework were studied by Sniatycki and Tulczjyew \cite{ST72} and 
Benenti \cite{Ben83, Ben11}. However, contributions from many other
authors play a central role for the development and understanding 
of generating functions and their relationship with the Hamilton-Jacobi 
theory, such as Viterbo \cite{Vit92}, Chaperon \cite{Cha95}
Maslov \cite{Mas65},  H\"ormander \cite{Hor71}, Weinstein \cite{Wei72} among many others.

In \cite{Jim15a}, the first author gives a strong argument for the construction of  
symplectic integrators
based on the fact that solutions of the Hamilton-Jacobi equation
for integrable and autonomous systems belong to a Lagrangian 
submanifold of the phase space \cite{AM78}. Starting from the 
classical approach, the product manifold of two copies of the phase space 
is created in \cite{Jim15a}. Then a generalized generating function with $4n$ variables 
is defined on the open ball of the product manifold centered at 
the initial condition (the source point of the map) and such that it 
contains the target point.
The generating function is directly associated to a primitive 1-form on the 
product manifold considered as a special symplectic manifold on 
the configuration space. The problem 
is translated from looking for the generating function, to looking for 
the 1-form, also known as the {\it Liouvillian form} \cite{Lib00}.
Two different Lagrangian
submanifolds arise, one defined by the generating function, solution of
the Hamilton-Jacobi equation, and the other invariant under the flow of the Hamiltonian vector
field. Applying an analogous argument to the Hamilton's method of characteristics, the
flow is searched in a transversal direction to the former submanifold. 
The suplementary space, transversal to the tangent space of the 
first Lagrangian submanifold, is projected by the induced projection\footnote{The 
induced projection is the one we used to define the symplectic form on the 
product manifold by its pull-back.}, onto a $2n$ subspace, where the 
original Hamiltonian system is finally evaluated. 
The projection induces a family of 
one step implicit symplectic integrators of generic order 1\footnote{For some 
particular values of the parameters, we observe an increment in the order 
of convergence (see Figure \ref{fig:6}).}, closely related to those already 
studied 
by Kang and co-workers and recently revisited by 
Xue and Zanna \cite{XZ14}. 

The main difference is that our numerical schemes contain some additional
terms which vanish when the stepsize goes to zero. In particular, the
symmetric integrators of the new family depend on $2n$ parameters and contain,
as a special case, the
midpoint rule. 
The goal of the present article is to perform a numerical study investigating 
those additional terms and their influence on the accuracy and performance.
It is a numerically oriented paper; for the geometrical
arguments and theoretical development we refer the reader to \cite{Jim15a,Jim15e}.
More theoretical results related to the numerical error and the claim by Ge and Marsen 
\cite{GM88,Ge91} about the impossibility of constructing exact numerical symplectic mappings 
is addressed in \cite{Jim15c}.
In fact, Ge's result \cite{Ge91} implies that our scheme should be exact up to numerical 
computer error and the residual error must be associated with the solution of the 
implicit equations.

\section{Hamiltonian systems and Symplectic integrators}
In what follows, we assume the reader is familiar with the terminology of differential geometry 
and vector bundles. For an introduction, the reader is referred to \cite{AM78,LM87,MR91}.

The approach presented in this work is based on the fact that 
Hamiltonian mechanics relies on the geometrical properties present
in the evolution of  
a mechanical system which accepts a Hamiltonian description. 
For this reason, in this section we use 
the standard notation in modern symplectic
geometry, topology and fiber bundles. 
The goal is to give a formal
geometrical framework to study symplectic integrators as 
isometries of a generic symplectic form on generic symplectic manifolds. 
In this way, the construction
of our symplectic methods takes a more abstract and general 
point of view solving, or correcting, some misunderstanding
arising when the analysis is restricted to symplectic vector 
spaces, or cotangent bundles of linear spaces. In general,
the starting point in symplectic integrator's analysis is the identification
of a cotangent bundle with a symplectic vector space by the isomorphisms 
\begin{eqnarray*}
    T^*(\mathbb R^n) \stackrel{iso}{=} (\mathbb R^n)^*\oplus\mathbb R^n \stackrel{iso}{=} \mathbb R^{2n},
\end{eqnarray*}
where $T^*(\mathbb R^n)$ is the cotangent bundle and 
$(\mathbb R^n)^*$ is the dual space of $\mathbb R^n$.
However, this point of view hides the
geometrical background of generating functions for constructing 
symplectic maps. Let us start with the main definitions and results.

A \emph{symplectic manifold} is a $2n$-dimensional manifold $M$ equipped with 
a non-degenerated, skew-symmetric, closed 2-form $\omega$, such that
at every point $m\in M$, the tangent space to $M$ at $m$, denoted
$T_mM$, has the structure of a symplectic vector space. One of the
basic properties in symplectic geometry is given by Darboux's
theorem which states that any symplectic manifold
is locally \emph{symplectomorphic} to a symplectic vector space $(V,\omega_0)$
with the canonical symplectic form $\omega_0(\cdot,\cdot)=\langle\cdot,J_0\cdot\rangle$,
where $\langle\cdot,\cdot\rangle$ is the canonical Euclidean structure
on $V$ and $J_0$ is the \emph{almost complex structure} represented by 
the matrix:
\begin{eqnarray*}
    J_0=\left(
    \begin{array}[h]{cc}
        0_n & -I_n\\
        I_n & 0_n
    \end{array}
    \right), \qquad 0_n,I_n\in M_{n\times n}(\mathbb R).
\end{eqnarray*}
$J_0$ is also known as the \emph{canonical symplectic matrix} on 
$V$.
Consequently, the tangent space to $M$ at $m$, with its symplectic 
form $\omega|_m=\omega_m$ and a suitable change of coordinates, can be completely described 
by the symplectic vector space with canonical symplectic coordinates
$({\bf  x},{\bf y})\in V$, $y_i=J_0x_i$. Darboux's theorem means that we systematically identify 
\begin{eqnarray*}
    (T_mM,\omega_m) \stackrel{iso}{=} (V,\omega_0), \qquad \forall m\in M,
\end{eqnarray*}
selecting the right change of symplectic coordinates to describe the 
dynamics on $M$ by the canonical symplectic coordinates $({\bf x},{\bf y})\in V$. 

\begin{remark}
Unfortunately, Darboux's theorem hides a very rich environment suitable
for investigating symplectic maps by use of generating functions and Liouvillian forms. In practice, it 
locally identifies all the symplectic manifolds of the same dimension with the cotangent
bundle. The information about the symplectomorphism which maps the symplectic 
manifold of interest to the cotangent bundle is lost when we apply Darboux's 
theorem. To avoid this loss of information we stay 
in the generic geometric framework of symplectic geometry by using special 
symplectic manifolds.
\end{remark}

In this geometrical framework, a Hamiltonian system $(M,\omega,X_H)$ is a vector field
$X=X_H$ on the symplectic manifold $(M,\omega)$ such that\footnote{Some 
authors write $i_{X_H}\omega = dH$ instead of 
(\ref{eqn:Ham:Def}), but it depends on the definition of $\omega$ as 
2-form and the choice of the complex structure $J$.} 
\begin{eqnarray}
   i_{X_H}\omega = -dH,
   \label{eqn:Ham:Def}
\end{eqnarray}
for a differentiable function $H:M\to\mathbb R$. 

A natural diffeomorphism $\flat:TM\to T^*M$ between the tangent
and the cotangent bundles of $M$ is given by the contraction of the symplectic 
form with the vector field in the following way
\begin{eqnarray*}
    X \mapsto \omega(X,\cdot). 
\end{eqnarray*}
The inverse of $\flat$ is denoted by 
$\sharp:T^*M\to TM$. Using $\sharp$, equation (\ref{eqn:Ham:Def}) is
written in vector field form as
\begin{eqnarray}
   X_H = J\nabla H,
   \label{eqn:Ham:Def2}
\end{eqnarray}
where  $\nabla$ is the standard gradient associated to the Euclidean structure.
Expression (\ref{eqn:Ham:Def2}) is best suited for applications. The 
equations of evolution can be written as
\begin{eqnarray}
    \dot \bz = J\nabla_{\bz} H(\bz),\qquad \bz\in M.
    \label{eqn:Field}
\end{eqnarray}

\begin{remark}
  When $M=T^*\mathbb R^n$ is equipped with a canonical symplectic
  basis in cotangent coordinates $\bz=(\bq,\bp)\in T^*\mathbb R^n$, the 
  Hamiltonian vector field is given by Hamilton's equations:
  \begin{eqnarray}
     \dot \bq = \frac{\partial H}{\partial \bp}(\bq,\bp),\qquad 
      \dot \bp = -\frac{\partial H}{\partial \bq}(\bq,\bp).
      \label{eqn:Ham}
  \end{eqnarray} 
\end{remark}

Poincar\'e discovered that the flow of any Hamiltonian vector 
field\footnote{Poincar\'e used the name of the ``fundamental problem of dynamics'',
which is to find the solutions of the ``fundamental equations of dynamics in
canonical form'' \cite{Poi99}.}
is a 1-parameter subgroup of symplectic diffeomorphisms.
Denoting such a flow by $\phi^t_H$, this implies that 
for each fixed $h\in\mathbb R$, $\phi^h_H$ is a symplectic map.

Let $\bz_0\in M$ be a point on the symplectic manifold and $\bz(t)$ the
integral curve to $X_H$ such that $\bz_0=\bz(0)$. By definition of the flow,
the mapping 
\begin{eqnarray*}
    \bz(t+h) = \phi^h_H (\bz(t))
\end{eqnarray*}
will propagate the solution from time $t$ to time $t+h$. 
A \emph{symplectic algorithm} with stepsize $h$, is the numerical 
approximation $\psi_h$ of the Hamiltonian
flow $\phi^h_H:M\to M$, which is an isometry of the symplectic form
$\omega$ . 
Specifically, consider the exact solution $\bz(t)$ of a Hamiltonian system
for the time $t\in [0,T]$, a discretization $\{t_i\}_{i=0}^N$ such that
$t_0=0$, $t_N=T$, $h=T/N=t_{n+1}-t_n$, and denote $\bz_n=z(t_n)$ for $0\le n \le N$. 
Let $U\subset M$ be a convex open neighbourhood of $\bz_n$ containing the target point
$\bz_{n+1}$. 

With these hypotheses, we define a \emph{symplectic integrator} 
as a map 
\begin{eqnarray*}
    \psi_h: U\subset M & \to & U\\
    \bz_n & \mapsto & \bz_{n+1} = \psi_h(\bz_n)
\end{eqnarray*}
smooth with respect to $h$ and $H$, and such that $\psi_h^*\omega = \omega$, where 
$\psi_h^*$ is the pullback of $\psi_h$ defined by 
\begin{eqnarray}
    (\psi_h^*\omega)_{\bz}(\bv,\bu) = \omega_{\psi(\bz)}(T\psi_h(\bv), T\psi_h(\bu)),\quad \bu,\bv\in T_zM.
    \label{eqn:def:pull}
\end{eqnarray}

\begin{remark}
    The vectors $T\psi_h(\bv)$ and $T\psi_h(\bu)$ belong to the tangent space 
    $T_{\psi_h(\bz)}M$ which, in general, is different from $T_{\bz}M$. Once we identify 
    $\omega_{\bz}$ and $\omega_{\psi_h(\bz)}$ with $\omega_0$ and the tangent spaces 
    $T_{\bz}M$ and $T_{\psi_h(\bz)}M$ with $V$, 
    condition (\ref{eqn:def:pull})
    becomes
    \begin{eqnarray*}
        \langle \bv,J \bu\rangle = \langle T\psi_h(\bv), J T\psi_h(\bu)\rangle.
    \end{eqnarray*}
    or equivalently 
    \begin{eqnarray*}
        \frac{\partial \psi_h}{\partial \bz}^T  J  \frac{\partial \psi_h}{\partial \bz}=J.
    \end{eqnarray*}
    which is the well-known symplecticity condition for $\mathbb R^{2n}$ viewed 
    as symplectic vector space. 
    It is worth nothing that in the last
    expression, $z$ are already local coordinates.
\end{remark}

In an analogous way, we define an \emph{implicit symplectic integrator} 
as a map 
\begin{eqnarray*}
    \varphi_h: U\times U & \to & U\\
    (\bz_n, \bz_{n+1}) & \mapsto & \bz_{n+1} = \varphi_h(\bz_n, \bz_{n+1})
\end{eqnarray*}
smooth with respect to $h$ and $H$, and such that $\varphi_h^*\omega = \omega$.

Since we are discretizing a flow, it is possible to consider an intermediate point $\bar{\bz}\in U$
 and two maps\footnote{In fact, they must be symplectic maps to have a consistent 
 symplectic integrator \cite{Jim15e}.} 
$\psi_1, \psi_2:U\to U$ such that $\bar{\bz}=\psi_1(\bz_n)=\psi_2(\bz_{n+1})$.
They let us rewrite the implicit scheme as $\varphi_h(\bz_n,\bz_{n+1}) = \psi_2^{-1}\circ \psi_1(\bz_n)$.

The pullback becomes 
$\varphi_h^*\omega = (\psi_2^{-1}\circ \psi_1)^*\omega= \psi_1^*\circ (\psi_2^{-1})^* \omega$ 
which produces the
corresponding symplecticity condition in the tangent spaces by
\begin{eqnarray}
    \left( \frac{\partial \psi_1^{-1}}{\partial \bar{\bz}}\right)^T  J 
    \left( \frac{\partial \psi_1^{-1}}{\partial \bar{\bz}}\right) =
    \left( \frac{\partial \psi_2^{-1}}{\partial \bar{\bz}}\right)^T  J  
    \left( \frac{\partial \psi_2^{-1}}{\partial \bar{\bz}}\right).
    \label{eqn:impl:symp}
\end{eqnarray}
Note that on the left hand side of (\ref{eqn:impl:symp}), $J$ is an endomorphism on $T_{\bz_n}U$ and 
on the right hand side, $J$ is an endomorphism on $T_{\bz_{n+1}}U$.

Condition (\ref{eqn:impl:symp}) says nothing about the mappings $\psi_1$ and $\psi_2$, but
only that the composition is a symplectic map. The reader interested
is referred to \cite{Jim15e} for a deeper discussion on conditions imposed on
$\psi_1$ and $\psi_2$.

The discrete scheme $\psi_h$ is said
of order $r\in\mathbb N$ if, as $h\to 0$ 
\begin{eqnarray*}
    \|\phi^h_H(\bz_n)-\psi_h(\bz_n)\| &=& \mathcal O(h^{r+1}), \qquad \bz_n=\bz(t_n)\in U.
\end{eqnarray*}

There are several methods for constructing symplectic integrators that reduce to finding symplectic maps 
$\phi^h_H$ between two different (closed) points on the integral curves of the vector field $X_H$. 
Here we are interested in the method of generating functions that we describe in the
next section.

\section{Generating functions}

It seems
that Poincar\'e was the first author who systematically studied the process of obtaining 
symplectic maps by generating functions. However, the terminology was different in his work:
Hamiltonian equations were called fundamental equations of dynamics in canonical form,
symplectic maps were called canonical transformations\footnote{Transformations which preserves
the canonical form of the fundamental equations of dynamics \cite{Poi99}.}. Generating functions
were the fundamental tool in his theory of integral invariants \cite{Poi99}.
Poincar\'e used these transformations to study periodic orbits of second genus 
in celestial mechanics. It may explain why this technique was not acknowledged immediatelly
by the numerical community. 
From the numerical point of view, generating functions for symplectic integrators
were systematically studied by Feng's team in the mid '80s
\cite{Kang85a, Kang85b, KG88, KHMD89} and later, among many others, 
by Ge and co-workers in late '80s and
'90s \cite{GM88,Ge90,Ge91,GD95}. Both the analytical and the numerical procedures  
use the same framework that we outline now.

Let $(M_1,\omega_1)$ and $(M_2,\omega_2)$ be two symplectic manifolds of the same 
dimension. A map $\phi:M_1\to M_2$ is 
called \emph{symplectic} if $\phi^*\omega_2=\omega_1$,
where, in general $\omega_1$ and $\omega_2$ are different symplectic forms.
In our case, $(M_1,\omega_1)=(M_2,\omega_2)$ are two copies of the same
symplectic manifold, but we will preserve the subindices to keep record
of what copy we refer to at every time.

Consider the product manifold $\tilde M=M_1\times M_2$ with canonical projections 
$\pi_i:\tilde M\to M_i$ for $i=1,2$, and define a two-form $\omega_{\ominus}$ on $\tilde M$ by 
\begin{eqnarray}
    \omega_{\ominus} &=& \pi_1^*\omega_1 - \pi_2^*\omega_2
    \label{eqn:def:sym}
\end{eqnarray}
We have the following results (see \cite[sec 5.2]{AM78} for the proofs):
\begin{itemize}
    \item $(\tilde M,\omega_{\ominus})$ is a symplectic manifold of dimension $4n$.
    \item for any symplectic map $\phi:M_1\to M_2$, the graph of $\phi$, denoted by $\Gamma_\phi$,
        and defined as
        \begin{eqnarray*}
            \Gamma_\phi&=& \left\{ \left(\bz,\phi(\bz)\right)\in \tilde M\ 
		|\ \bz\in M_1, \phi(\bz)\in M_2 \right\},
        \end{eqnarray*}
        is a Lagrangian submanifold of $\tilde M$. This means $\omega_{\ominus}|_{\Gamma_\phi}=0$
    \item Since $\theta_\ominus$  is a closed 1-form by the identity 
        $d\theta_{\ominus}|_{\Gamma_\phi}\equiv 0$, using Poincar\'e's lemma,
        $\theta_{\ominus}$ is also an exact 1-form on $\Gamma_\phi$. Then, 
        there exists a function $S$ defined on
        the Lagrangian submanifold $\Gamma_{\phi}$ such that its differential concides with the
        restriction of the 1-form $\theta_{\ominus}$ to $\Gamma_{\phi}$
        \begin{eqnarray}
            S:\Gamma_\phi\to \mathbb R, \qquad{\rm such\ that}\quad dS|_{\Gamma_\phi}\equiv \theta_{\ominus}|_{\Gamma_\phi}.
            \label{eqn:dif:func}
        \end{eqnarray}
	$S:\Gamma_{\phi}\to\mathbb R$ is called a \emph{generating function} for the symplectic
	map $\phi$.
\end{itemize}

Symplectic maps, generating functions and Lagrangian submanifolds 
are closely related. For instance,  in a generic symplectic manifold $M$, any Lagrangian 
submanifold $\Lambda\subset M$ 
which is transverse to 
the fibers of the projection $\pi_{M}:T^* M\to M$ 
can be locally parameterized by a suitable atlas of 
(local) functions $S_i$ for short times \cite{Hor71, Mas65, Car22}. 
Since symplectic integrators are mappings close to 
the identity ($h$ small), we are not concerned with the global behaviour 
of the Lagrangian submanifolds and all our analysis will be local.

The standard procedure of the method of generating functions \cite{Kang85a, Kang85b} is as follows:
1) look for a suitable 1-form $\theta\in T^*\tilde M$ such that 
$d\theta = \omega_{\ominus}$;
2) obtain the Lagrangian submanifold $\Lambda\subset \tilde M$, associated with $\theta$
and a function $S:\Lambda\to\mathbb R$ satisfying $dS=\theta$;
3) solve the Hamilton-Jacobi equation on the Lagrangian submanifold;
4) design a numerical method for approximating such a solution giving the mapping 
$\bz_{n+1}=\psi_h(\bz_n)$ for $\bz_n\in M_1$ and $\bz_{n+1}\in M_2$.

Kang's procedure for solving the Hamilton-Jacobi equation
is to approximate the generating function using Taylor series expansions \cite{Kang85a}. 
Menyuk  used the Picard iteration to obtain such a function \cite{Men84} and 
Channell and Scovel used symbolic computing software \cite{CS90}.
As pointed out by Miesbach and Pesch in \cite{MP92}, these approaches require 
higher order derivatives of the Hamiltonian which complicates the final scheme.

To keep the notation simple, from now on lowercase variables belong to $M_1$
and uppercase ones belong to $M_2$, in particular
$\bz = \bz_n$ and $\bZ=\bz_{n+1}$. We use also Einstein notation in the 
definition of differential forms, i.e., $p_idq_i:=\sum_{i=0}^n p_idp_i$.

\subsection{The alternative method of Liouvillian forms}
Recently, the first author made some contributions on the subject \cite{Jim15a},
where he avoids the solution of the Hamilton-Jacobi equation and recovers 
a numerical algorithm from geometrical interpretation of the solutions. 
Indeed, solutions of the Hamilton-Jacobi 
equation on the product manifold $\tilde M$, belong to a Lagrangian submanifold. 
This submanifold can be related to the characteristic bundle of the Hamiltonian
vector field generated by a generalized generating function on $\tilde M$.
The differential of this generalized function is actually a Liouvillian form,
\emph{i.e.} a 1-form $\lambda$ on $\tilde M$ such that $d\lambda=\omega_\ominus$. 
A suitable projection $\pi:\tilde M\to N$ to a $2n$ dimensional submanifold 
$N\subset \tilde M$ gives the right point where we must evaluate the discrete
flow in order to have a symplectic integrator. Contrasting this construction with 
Kang's procedure, we look for good values of the discrete flow on the Lagrangian 
submanifold defined by the Liouvillian form before the projection, instead of 
approximating the projection of the solution by Taylor series.
Consequently, points 3) and 4) are not relevant in this construction. 
For more details the reader is referred to  \cite{Jim15a}.

Our point of departure is the 
family of Liouvillian forms $\theta_{\alpha,\beta,\gamma}$ constructed in \cite{Jim15a}. 
Consider two copies of the phase space $M_1=T^*\mathcal Q_1$, $M_2=T^*\mathcal Q_2$ and define 
the product manifold $\tilde M= M_1\times M_2$ equipped with the 2-form $\omega_\ominus$ 
given in (\ref{eqn:def:sym}). 
As we know, $(\tilde M,\omega_\ominus)$ is a symplectic manifold of dimension $4n$. Define 
local coordinates $(q_i,p_i)\in M_1$ and $(Q_i,P_i)\in M_2$, $i=1,\cdots,n$, in each copy 
of the phase space and consider real numbers $\alpha,\beta,\gamma\in\mathbb R$. With this notation, the family of 
primitive 1-forms $\theta_{\alpha,\beta,\gamma}$ is given by
\begin{eqnarray}
    \theta_{\alpha,\beta,\gamma} &=& \alpha(p_i dq_i + Q_idP_i) - (1-\alpha)(P_idQ_i - q_idp_i)  \nonumber \\
                & & + \beta(q_idq_i + Q_idQ_i) - \gamma(p_idp_i +  P_i dP_i),
                \label{eqn:theta1}
\end{eqnarray}
Since $\omega_\ominus = d\theta_{\alpha,\beta,\gamma}$ for any values of the parameters,
we can give a full set of 3$n$ different parameters $\{\alpha_i,\beta_i,\gamma_i\}_{i=0}^n$  
for every 1-form $\theta_{\alpha,\beta,\gamma}$\footnote{The case $\beta_i=\gamma_i=\sqrt{\alpha_i(1-\alpha_i)}$ with $\alpha_i\in[0,1]$ corresponds to 
the family $\theta_{\phi_i}$ constructed in \cite{Jim15a} from the more elementary symplectic rotation 
$\alpha_i=\cos^2(\phi_i)$, $i=1\cdots,n$, on $\tilde M=M_1\times M_2$.}.
Note that the elements associated with $\beta$ and $\gamma$ belong to the 
kernel of the differential and they are known as the \emph{gauge} elements of 
$\theta_{\alpha,\beta,\gamma}$. They do not modify the symplectic structure of
$(\tilde M, d\theta_{\alpha,\beta,\gamma})$ but the solution in the Lagrangian 
submanifold will be, in general, different for each combination of parameters.

\begin{remark}
Every point $(\alpha_i,\beta_i,\gamma_i)\in \mathbb R^{3n}$, $i=0,\cdots,n$,
is associated to a 1-form whose differential is exactly the symplectic form $\omega_\ominus$.
However, good values for a symplectic map approximating the
Hamiltonian flow in a suitable way, belong to an open ball around the 
point $(\alpha_i,\beta_i,\gamma_i)=(1/2,0,0)$, associated to the mid-point 
symplectic map. 
In Section \ref{sec:num} we will find values 
for these parameters for a concrete example.
\end{remark}

All members of the family $\theta_{\alpha,\beta,\gamma}$ are 1-forms, locally closed on the graph
$\Gamma_\phi\subset \tilde M$ of a generic symplectic map $\phi:(M_1,\omega_1)\to (M_2,\omega_2)$. 
By Poincar\'e's lemma,
there exists a generating function $S=S_{\alpha,\beta,\gamma}$ depending on
$(\alpha,\beta,\gamma)$ such that
$dS=\theta_{\alpha,\beta,\gamma}$.
The Lagrangian submanifold $\Lambda_{\alpha,\beta,\gamma}$ parameterized in local coordinates by the 
equation $dS=\theta_{\alpha,\beta,\gamma}$ has explicit form 
\begin{eqnarray}
   \frac{\partial S}{\partial q_i} &=& \alpha p_i + \beta q_i \\ 
   \frac{\partial S}{\partial p_i} &=& - (1-\alpha) q_i - \gamma p_i \\
   \frac{\partial S}{\partial Q_i} &=& -(1-\alpha) P_i + \beta Q_i \\
   \frac{\partial S}{\partial P_i} &=& \alpha Q_i - \gamma P_i, 
\end{eqnarray}
and must satisfy the homogeneous Hamilton-Jacobi equation 
$H\left(\Lambda_{\alpha,\beta,\gamma}\right) = 0$.

Consider a 2n-dimensional submanifold of the product manifold $N_{\alpha,\beta,\gamma}\subset \tilde M$ 
with coordinates $(\bar Q_i, \bar P_i)$, and define the projection $\pi_S:\tilde M\to N_{\alpha,\beta,\gamma}$ 
by 
\begin{eqnarray}
  \pi_S = J\circ(\pi_1-\pi_2)(\Lambda_{\alpha,\beta,\gamma}).
  \label{eqn:proj:1}
\end{eqnarray}

The submanifold $N_{\alpha,\beta,\gamma}$ is given in local coordinates by the equations
\begin{eqnarray}
    \bar{Q_i} &=&  \alpha Q_i + (1-\alpha) q_i + \gamma( p_i - P_i) \nonumber\\
    \bar{P_i} &=&  \alpha p_i + (1-\alpha)P_i + \beta( q_i - Q_i ).
    \label{eqn:proj}
\end{eqnarray}
Note that on the diagonal $\Delta_{\tilde M}$ of $\tilde M$ defined by
\begin{eqnarray*}
    \Delta_{\tilde M}=\left\{(q_i,p_i,Q_i,P_i)\in\tilde M\ |\ q_i=Q_i,\ p_i=P_i,\forall i=1\cdots n\right\},
\end{eqnarray*}
the submanifold $N_{\alpha,\beta,\gamma}$ is a standard symplectic submanifold with the induced 
symplectic form $\omega_\Delta$ which fulfils $\omega_\ominus=\pi^*_S\omega_{\Delta}$.
Moreover, given the inclusion 
\begin{eqnarray}
  i:M & \to &  \tilde M\\
  (q_i,p_i)&\mapsto&(q_i,p_i,q_i,p_i)
\end{eqnarray}
whose image is exactly the diagonal $\Delta_{\tilde M} = i(M,\omega)$
we have $\omega=i^*\omega_{\Delta}$. 

\subsection{The implicit symplectic integrators}
For $(\bz,\bZ)\in\tilde M$ where $\bz=(q_i,p_i)\in M_1$ and $\bZ=(Q_i,P_i)\in M_2$,
the family of symplectic integrators given in \cite{Jim15a} is obtained by 
the following implicit scheme\footnote{Note that the projection (\ref{eqn:proj:1}) with
local coordinates (\ref{eqn:proj}) is linear
in $(q_i,p_i,Q_i,P_i)$ and consequently it coincides with its linearization $T\pi_S$.}
\begin{eqnarray}
    \bZ &=& \bz + h J\nabla H( T\pi_S(\bz,\bZ) ) 
    \label{eqn:sympl}
\end{eqnarray}
or in extended form
\begin{eqnarray}
    Q_i &=& q_i + h  \frac{\partial H}{\partial p_i}\left(
        \alpha Q_i + \alpha^\prime q_i + \gamma( p_i - P_i),
        \alpha p_i + \alpha^\prime P_i + \beta( q_i - Q_i )
    \right)\\ 
    P_i &=& p_i - h  \frac{\partial H}{\partial q_i}\left( 
        \alpha Q_i + \alpha^\prime q_i + \gamma( p_i - P_i),
        \alpha p_i + \alpha^\prime P_i + \beta( q_i - Q_i )
    \right) \nonumber 
    \label{eqn:sympl:pq}
\end{eqnarray}
where we used $\alpha^\prime=1-\alpha$ to simplify the expressions.
We shall return to these expressions later, when we study the symmetric
case. At this stage, it is worth to note that for 
$(\alpha,\beta,\gamma)=(0,0,0)$
and $(\alpha,\beta,\gamma)=(1,0,0)$, we obtain the staggered Euler symplectic 
integrators A and B \cite{HLW10}
\begin{eqnarray}
    Q_i = q_i + h  \frac{\partial H}{\partial p_i}\left( q_i, P_i \right) &\quad&
    Q_i = q_i + h  \frac{\partial H}{\partial p_i}\left( Q_i, p_i \right) \\
    P_i = p_i - h  \frac{\partial H}{\partial q_i}\left( q_i, P_i \right) &\quad&
    P_i = p_i - h  \frac{\partial H}{\partial q_i}\left( Q_i, p_i \right) 
    \label{eqn:sympl:pqStag}
\end{eqnarray}

The main remark on the choice of the method (\ref{eqn:sympl}) is that projection $\pi_S$ given in (\ref{eqn:proj:1}) is a 
``rearrangement'' of the data coded in the 
Lagrangian submanifold $\Lambda_{\alpha,\beta,\gamma}$ in a close to symplectic submanifold 
$N_{\alpha,\beta,\gamma}$. 
This choice is rooted in the notion of geometrical solution to the 
Hamilton-Jacobi equation 
\begin{eqnarray*}
    \frac{\partial S}{\partial t} + H(\lambda) = 0
\end{eqnarray*}
where the function $H:M\to\mathbb R$ of the Hamilton-Jacobi equation is evaluated on 
a Lagrangian submanifold $\lambda\subset M$ of a symplectic 
manifold $(M,\omega)$ given by local equations 
\begin{eqnarray*}
    \lambda=\left\{ (\bq,\bp)\in M\ |\ \bq=\bq,\ \bp=\frac{\partial S}{\partial \bq} \right\},
\end{eqnarray*}
and differs from the energy Hamiltonian
function, only by a constant (normally denoted by $-E$). For the non-evolutionary 
Hamilton-Jacobi equation, considering the characteristic 
bundle, we evaluate the Hamiltonian vector field on the sub-bundle $(T\lambda)^\perp\subset TM$
transversal to $T\lambda$.
An example is provided by the staggered Euler integrator of type A, with the 
difference approximation 
\begin{eqnarray}
    \bZ &=& \bz + h J\nabla H( \bQ, \bp ),
\end{eqnarray}
of the differential equation $\dot \bz = J\nabla H(\bQ,\bp)$ associated to  
the Hamilton Jacobi equation $H(\bq,\bP)=0$
for $\bz=(\bq,\bp)$ and $\bZ=(\bQ,\bP)$.

\section{Numerical examples}
All the  members of the family (\ref{eqn:sympl:pq}) for which 
$\alpha\in (0,1)$ produce implicit 
symplectic integrators that we can compute in an iterative predictor-corrector
scheme. 

For the tests, we consider the perturbed
Hamiltonian pendulum 
as toy example.
The algorithm was implemented in python for testing
the accuracy of the method and compared with
the explicit method $\mathcal {SABA}_2$ with coefficients 
$$\{a_0=a_2=1/6, b_0=b_1=1/2, a_1=2/3\},$$
studied by Laskar and Robutel in \cite{LR01}. 
We used the staggered symplectic Euler integrator as 
predictor and we compute, iteratively, the intermediate coordinates $\bar \bz=(\bar \bQ, \bar \bP)$
which are used to compute the value of $\bz_{n+1}=(\bp_{n+1}, \bq_{n+1})$ in $\kappa$
iterations.

The general structure of the algorithm is the following 
\begin{center}
  \begin{tabular}{|rl|}
    \hline
    \multicolumn{2}{|l|}{
        {\bf Algorithm 1.} }\\
    \hline
        & $!$Setup the predictor\\
     1: & $\bz^{[0]}=\bz_n + hJ\nabla H(\bz_n)$\\
     2: & for $j=0:\kappa$ do\\
        & $\quad$ $!$compute the point $\bar \bz=(\bar \bQ,\bar \bP)$ \\
     3: & $\quad \bar{\bQ} =  \alpha \bq^{[j]} + (1-\alpha)\bq_n + \gamma( \bp_n - \bp^{[j]})$ \\
     4: & $\quad \bar{\bP} =  \alpha \bp_n + (1-\alpha)\bp^{[j]} + \beta( \bq_n - \bq^{[j]} ).$ \\
        & $\quad$ $!$compute the corrector\\
     5: & $\quad$ $\bz^{[j+1]}= \bz_n + hJ\nabla H(\bar \bz)$ \\
     6: & end for\\
     7: & $\bz_{n+1}=\bz^{[\kappa]}$\\
    \hline
  \end{tabular}
\end{center}

\begin{figure}[h!]
 \centering
 \includegraphics[width=5in]{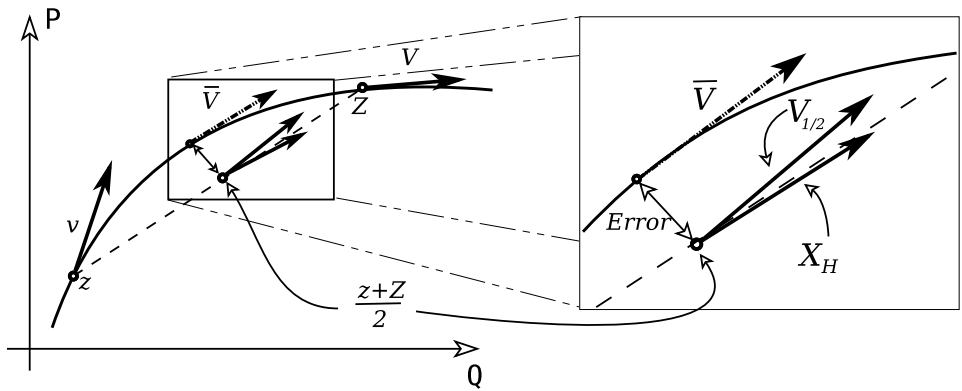}
 \caption{A simple illustration of the natural error obtained by the 
 symplectic mid-point rule. In this plot we assume known the vectors 
 $v=X_H(\bz)=\dot \bz$ and $V=X_H(\bZ)=\dot \bZ$ and we try to approximate 
 $\bar V$ by computing $X_H\left( \frac{\bz+\bZ}{2}\right)$. 
 Note that it coincides with $V_{\frac12}=\frac12(v+V)$ only if
 $X_H$ is a linear vector field. Moreover, it coincides with $\bar V$ only 
 when the solution $\bz(t)$
 is linear in $t$, or equivalently when $X_H$ is constant.  }
 \label{fig:midpoint}
\end{figure}

\subsection{The simple Hamiltonian pendulum}
The simple Hamiltonian pendulum obeys the Hamiltonian equation
\begin{eqnarray}
  H(p,q) &=& \frac{1}{2}p^2 - \epsilon \cos{q}
  \label{eqn:pend}
\end{eqnarray}
where we have fixed $\epsilon=0.03$.
We performed several tests
for different values of $(\alpha,\beta,\gamma)$ and the step size $h$. The total number of 
steps $N$ used in the simulations was variable but almost all the results used $N=10^{5}$. 

A full set of generic parameters $(\alpha, \beta,\gamma)$ for our symplectic integrators in a
Hamiltonian system with $n$ degrees of freedom, has $3n$ elements. We are interested in some particular 
cases and the dimension of the space of parameters is different for each one. Since it is not
evident in a 1 degree of freedom problem, we enumerate the dimension of the space of parameters 
for the cases used in the tests: 1) the case $(\alpha=0.5,\beta=0,\gamma=0)$
has null dimension (is a point);
2) the case $(\alpha,\beta,\gamma=\beta)$ has dimension $2n$;
3) the case $(\alpha, \beta=\sqrt{\alpha(1-\alpha)},\gamma=\sqrt{\alpha(1-\alpha)})$ has dimension $n$;
4) the case $(\alpha=0.5, \beta=\sqrt{\alpha(1-\alpha)},\gamma=\sqrt{\alpha(1-\alpha)})$
has null dimension, it is the point $(0.5,0.5,0.5)$; 
5) the case $(\alpha=0.5,\beta,\gamma)$ has dimension $2n$. Below the reader will see that
we look for optimal values of $\beta$ and $\gamma$ in an open ball centered 
at the origin of $\mathbb R^{2n}$.

\subsubsection{Preliminary tests}
   We consider the classical symplectic mid-point rule as starting point since it is the 
   algorithm corresponding to $(\alpha,\beta,\gamma)=(0.5,0.,0.)$. The error associated to this
   case is easy to understand if we consider the projection of the orbit on the 
   phase space (Figure \ref{fig:midpoint}). What we are looking for is 
   some point $\bar \bz$, as close as possible to the real orbit, such that on it, 
   the real vector field is parallel to the numerical value. Of course,
   the point $\bar \bz$ coincides with the mid point $\frac{1}{2}(\bz_n + \bz_{n+1})$
   when the Hamiltonian vector field is constant.
   
\begin{figure}[h!]
 \centering
 \includegraphics[scale=0.54]{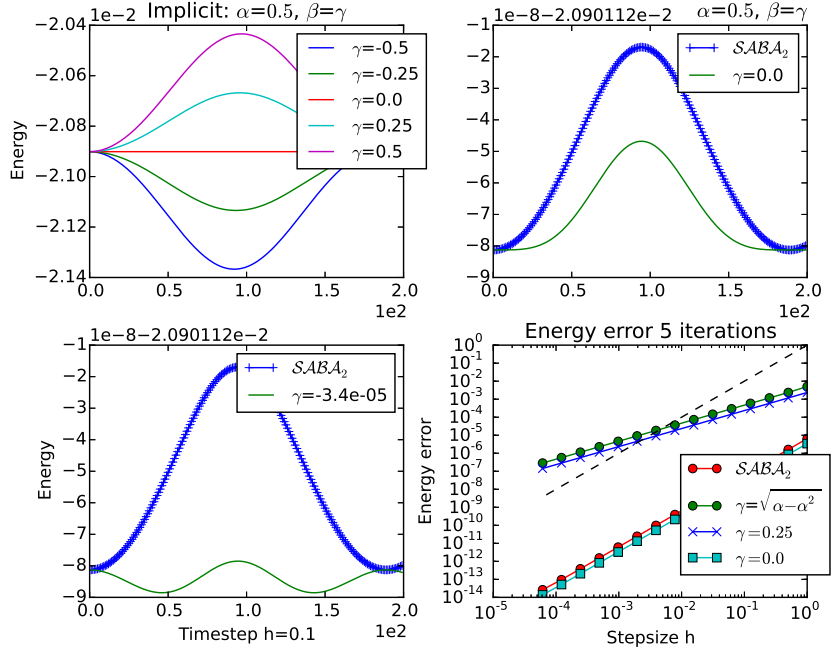}
 \caption{Comparison of the oscillations of the numerical solution around the 
 constant energy integral, between the new integrator with $\alpha=0.5$ 
 and the explicit $\mathcal{SABA}_2$ integrator. All plots have initial condition 
 $(q_0,p_0)=(0.8, 0.0)$, the stepsize in the three plot of energy integrals is $h=0.2$.
 Above-left: variation of the $\beta=\gamma$ parameters in $\{-0.5,-0.25,0, 0.25,0.5\}$
 Above-right: new method for $\beta=\gamma=0$ wompared to $\mathcal{SABA}_2$.
 Below-left: minimal oscillation around the energy integral of the new method for $\beta=\gamma=-3.4\times 10^{-5}$.
 Below-right: generic members of the family have order 1. The dotted line corresponds to $h^2$. }
 \label{fig:2}
\end{figure}

   When computing the symplectic mid-point rule with five iterations in the resolution of the 
   implicit scheme, and comparing with the explicit integrator $\mathcal {SABA}_2$ 
   from \cite{LR01}, 
   we found that the former is more accurate for the selected 
   set of initial conditions (upper-right plot in Figure \ref{fig:2}). 
   However, it is significantly more costly because 
   it is an iterative scheme.   

\begin{figure}[h!]
 \centering
 \includegraphics[width=5in]{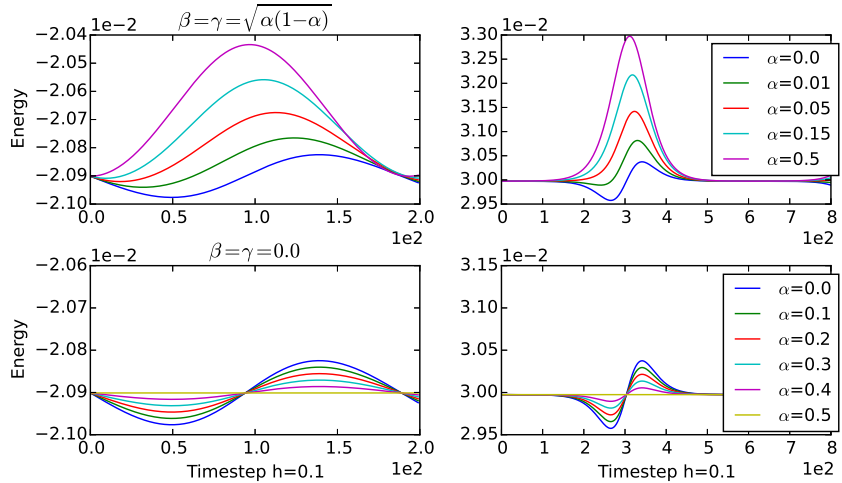}
 \caption{Oscillations around the energy integral of the numerical solutions 
 computed with the new implicit method. Left column corresponds to initial
 condition $(q_0,p_0)=(0.8,0.0)$ and $(q_0,p_0)=(3.1,0.0)$ for the right column.
 Different graphs in each panel correspond to the solution varying from 
 $\alpha=0.0$ to $\alpha=0.5$.  Above: progresive deformation of the 
 oscillations for the rotation $\beta=\gamma=\sqrt{\alpha(1-\alpha)}$.
 Below: the case $\beta=\gamma=0.0$.}
 \label{fig:ene1}
\end{figure}

With this error reference's framework, we consider the case where
$(\alpha,\beta,\gamma)$ are given by the simple symplectic rotation in \cite{Jim15a}.
The rotation corresponds to the parameters $\beta=\gamma=\sqrt{\alpha(1-\alpha)}$ 
for $\alpha\in[0,1]$. The relevant observations follow:

{\it a)} Variations in $\alpha$ from 0 to 1 result in a progressive deformation
of the orbit which connects consistently the symplectic Euler methods $A$ with $B$.
The energy is well conserved with the well-known oscillations for Euler methods.
The oscillations  
are of the same order of magnitude for every member of the family.
In the symmetric case $(\alpha,\beta,\gamma)=(0.5,0.5,0.5)$, the oscillations are really
large compared with the mid-point case $(\alpha,\beta,\gamma)=(0.5,0.,0.)$ 
(Figure \ref{fig:ene1}).

\begin{figure}[h!]
 \centering
 \includegraphics[scale=0.55]{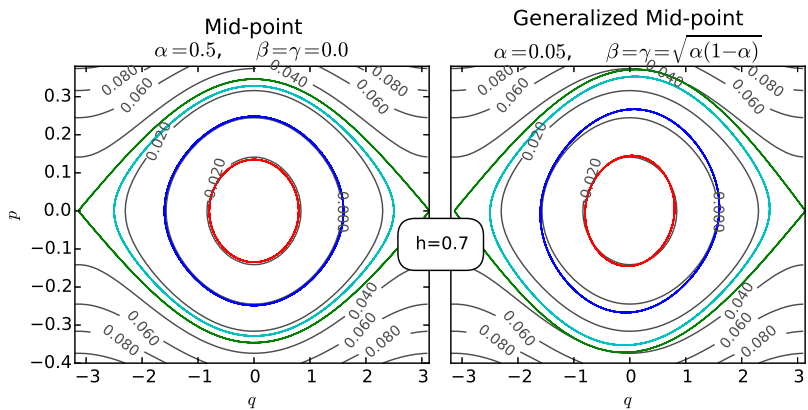}
 \caption{Orbits of the numerical solution for the perturbed Hamiltonian pendulum. We compare 
 two different members of the family of implicit integrators $(\alpha,\beta,\gamma)=(0.5,0,0)$
 and $(\alpha,\beta,\gamma)=(0.1,0.3,0.3)$ vs. the exact solution.
 There are four orbits in each plot, with initial conditions 
 $(q_0,p_0)\in\{(0.8,0.),(1.6,0.),(2.5,0.),(3.13,0.)\}$, 
 stepsize $h=0.7$ and $k=5$. Note the deformation of the orbits at the 
 right panel.}
 \label{fig:1}
\end{figure}

{\it b)} Tests for stepsizes $h$ between 0.001 and 0.1 exhibit a good behaviour 
and some deformation in the orbits is detected for $0.2<h<0.8$ and 
$\beta=\gamma=\sqrt{\alpha(1-\alpha)}$ (see right panel in 
Figure \ref{fig:1}). For $h\ge 1.0$ the orbits 
experienced high oscillations and for initial conditions close to the hyperbolic
fixed points the numerical solution went to the unbounded region.
Deformation of the orbits in phase space corresponds to the oscillation of the 
numerical solution of the energy integral 
around the exact constant value of $H$. On the other hand, for small values of $\beta$ and
$\gamma$ the numerical solutions has a good behaviour (see left panel in Figure \ref{fig:1})
and some tests for $h=10$ are really satisfactory.

The main remark from these preliminary tests is that parameters $\beta$ and $\gamma$ 
let us control the numerical 
error in the energy integral, and that values of the parameters close to 
$(\alpha,\beta,\gamma)=(0.5,0,0)$ have very small error. In the rest of the 
discussion we fixed $\alpha=0.5$ and we considered
small values for the other parameters. 
This value for $\alpha$ is strategic since $\alpha$ controls the symmetry of oscillations around the constant energy
for positive and negative directions in time. This is evident from the well-known fact that 
symplectic methods are reversible in time if and only if they are symmetric.
Moreover, oscillations for $\alpha\in(0,\frac12)$ when   
$\beta=\gamma=\sqrt{\alpha(1-\alpha)}$) in positive time $(\alpha,t)$, are symmetric 
with respect to those for $(1-\alpha,-t)$ in negative time.

\subsubsection{The symmetric case $\alpha=0.5$\label{sec:num}}
We perform this analysis in two separate cases: the first one when $\gamma$
is free and $\beta=\gamma$, and the second one when both $\beta$ and $\gamma$ are 
independent.

{\it a)} Let us fix $\alpha=0.5$ and $\beta=\gamma$, with $\gamma\in[-0.5,0.5]$. 
In this case, the amplitude of the oscillations in the energy 
goes from a positive to a negative phase (upper-left panel in Figure \ref{fig:2}).
Moreover, for $\gamma\sim 0$ and five iterations (for solving the implicit scheme) the
amplitude of the oscillations is lower than oscillations of $\mathcal{SABA}_2$
with a better behavior (lower-left panel in Figure \ref{fig:2}). This fact 
is remarkable since our implicit symplectic integrator is a one step method. Moreover
the tests on the behaviour of the error with respect to the stepsize
show a better accuracy than $\mathcal{SABA}_2$ (lower-right plot in Figure \ref{fig:2}).

A finer analysis for different orbits with initial conditions
close to one of the hyperbolic fixed points 
gives us new insight for the understanding of our numerical scheme.

\begin{figure}[h!]
 \centering
 \includegraphics[width=5in]{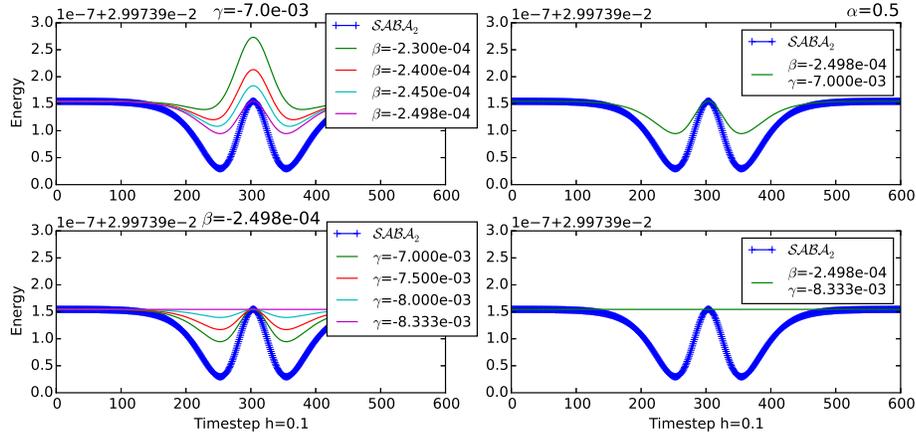}
 \caption{Search of optimal values for $\beta$ and $\gamma$ with $\alpha=0.5$. 
 Above: variation of the $\beta$ parameter changes the position of the
 central cusp. 
 Below-left: variation of the $\gamma$ parameter changes the position of the
 sourounding cusps. 
 Below-right: an optimal combination of parameters $\beta$ and $\gamma$ reducing the error 
 oscillations around $err<10^{-13}$. See left panel in Figure \ref{fig:4} for a zoomed
 version of this plot. 
  The fixed oscillating orbit corresponds
 to  the numerical solution of the $\mathcal{SABA}_2$ integrator. 
 All plots have initial condition 
 $(q_0,p_0)=(3.1, 0.0)$ and stepsize $h=0.1$.}
 \label{fig:3}
\end{figure}

{\it b)} Fix $\alpha=0.5$ and consider independent values for $\beta,\gamma\in[-0.5,0.5]$.
We observe in this case,
that each one of the parameters controls a different part of the oscillation around the 
energy integral. To show that, 
we choose an initial condition close to one hyperbolic fixed point. In an oscillation of a numerical 
solution given by a symmetric symplectic method (for instance $\mathcal{SABA}_2$), 
three cusps arise\footnote{We need to check if other methods and other 
Hamiltonian problems behave in a similar way.}. In our results, $\beta$ controls
the central cusp in a very clear and definite way. This property is shown 
in the upper row of Figure \ref{fig:3}. In this plot we fixed $\gamma=-00.7$ 
and we modified $\beta$ visually such that the cusp be close to the 
constant energy. Then, we fix $\beta$ and modify $\gamma$, which controls the two cusps
until we arrive at a very flat solution as we can see in the lower panel in
Figure \ref{fig:3}.

\begin{figure}[h!]
 \centering
 \includegraphics[width=5in]{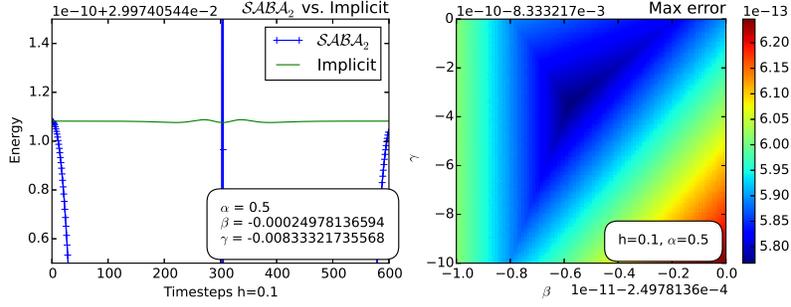}
 \caption{Search of the optimal values for $\beta$ and $\gamma$.  Left panel shows the oscillations
 of the error in the energy integral computed with values for  $\beta$, $\gamma$ found 
 by the plot in the right panel. Vertical curves correspond to the numerical 
 solution computed with $\mathcal{SABA}_2$. This is a zoomed version of the lower-right panel 
 in Figure \ref{fig:3}. Right panel shows 
 the distribution of the maximal variation of the error in a small region of the 
 $\beta-\gamma$ plane. The behaviour is equivalent for bigger and smaller regions:
 the minimal variation converges to a very particular region with radius around
 $10^{-13}$. }
 \label{fig:4}
\end{figure}

By a continuity argument we search for optimal values of $\beta$ and $\gamma$,
reducing the error in the numerical solution. This is done as a manual and visual 
process, computing the maximum variation of the error within an open square
around $(\beta,\gamma)=(0,0)$. We obtain a well defined region where 
such values belong. Maximum variation in the energy integral for the perturbed
Hamiltonian pendulum, with initial conditions $(q_0,p_0)=(3.1,0.0)$, 
stepsize $h=0.1$, perturbing term $\epsilon=0.03$ and $k=5$ iterations 
went to $err< 5\times 10^{-14}$, which is remarkable for a one step integrator. 
Figure \ref{fig:4} shows the maximum variation of the error for a subset of
$\mathbb R^2$ and the computed numerical energy compared with
$\mathcal{SABA}_2$. The vertical line in the middle of the left panel in 
Figure \ref{fig:4} is the central cusp of the $\mathcal{SABA}_2$ integrator. It is 
a zoomed version of the lower-right panel from Figure \ref{fig:3}.

\subsubsection{Interpretation of the parameters $\beta$ and $\gamma$} 
The fact that the new parameters $\beta$ and $\gamma$ control the numerical error 
is an important result. However, the way we obtained these values is quite heuristic.
In order for these parameters to be really advantageous, 
their behaviour with respect to the variations in the stepsize $h$ must be investigated. 

\begin{figure}[h!]
 \centering
 \includegraphics[width=5in]{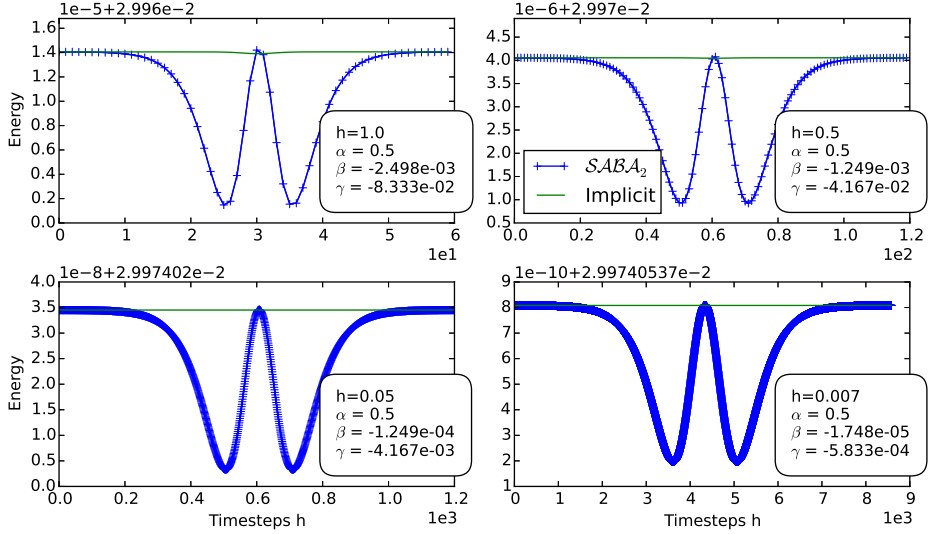}
 \caption{Energy integral of the numerical solutions of the new integrator compared to 
 the explicit $\mathcal{SABA}_2$ integrator with initial conditions $(q_0,p_0)=(3.1, 0.0)$
 and several stepsizes. For each value in the stepsize $h$, the values in the parameters 
 $\beta$ and $\gamma$ which minimizes the oscillations in the energy integral are proportional 
 to  $h$. Values of those parameters for the down-right panel were extrapolated by the 
 constants of proportionality obtained from the other panels. }
 \label{fig:5}
\end{figure}

\begin{figure}[h!]
 \centering
 \includegraphics[scale=0.55]{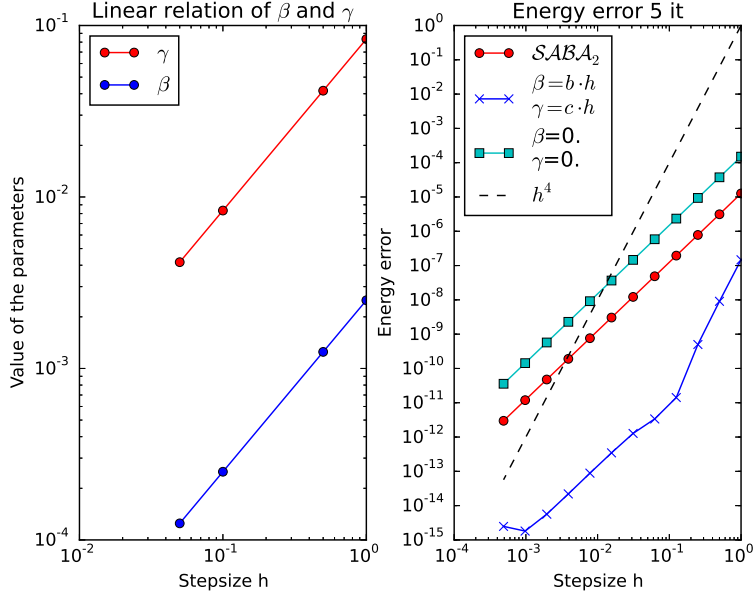}
 \caption{Left panel: The linear behaviour of the optimal values for
 $\beta$ and $\gamma$ as a function 
 of the stepsize $h$. Right panel, the optimal values of 
 $\beta$ and $\gamma$ for each $h$ increase the order of the integrator up to order four.}
 \label{fig:6}
\end{figure}

To check their dependency on the stepsize $h$, we performed several tests with the same 
initial conditions as in the previous tests. We searched for the optimal values 
of $\beta$ and $\gamma$ for $h\in\left\{ 0.05, 0.1, 0.5, 1.0 \right\}$ and we obtained an almost
perfect linear relationship $\beta=h\cdot b$ and $\gamma=h\cdot c$ for $b=-2.4978136594\times 10^{-3}$
and $c=-8.33321735568\times 10 ^{-2}$ (see upper row and lower-right panel in Figure 
\ref{fig:5}). To verify the linear dependence, we extrapolated and interpolated several values 
of the parameters for different $h$. The linear relationship holds 
in  a very accurate way. The lower-right panel in Figure \ref{fig:5} is computed
with the extrapolated parameters for $h=0.007$. 
Figure \ref{fig:6} shows the values found by visual inspection
for $h\in\left\{ 0.05, 0.1, 0.5, 1.0 \right\}$ and the line joining them.

This fact reveals an important property giving  a clue on the way a 
numerical symplectic integrator can preserve the energy (and any other) integral  
 in a very accurate way.
The extra elements in $(\bar Q,\bar P)$ from (\ref{eqn:proj}) can be written 
as an approximation of the gradient vector field $\nabla_\bz H(\bz)=-J\dot \bz$ 
considering 
\begin{eqnarray*}
    \gamma (p_i-P_i) &=&  -h\gamma \left(\frac{P_i-p_i}{h}\right)\\
    \beta(q_i-Q_i)  &=& -h\beta \left(\frac{Q_i-q_i}{h}\right)
\end{eqnarray*}
as the approximation $hA(-J\frac{\bZ-\bz}{h})$ where 
\begin{eqnarray*}
    A=\left(
    \begin{array}[h]{cc}
        \gamma I_n & 0_n\\
        0_n & -\beta I_n
    \end{array}
    \right), \qquad 0_n,I_n\in M_{n\times n}(\mathbb R).
\end{eqnarray*}

With this notation, the symplectic map (\ref{eqn:sympl}) is redefined  
in terms of the numerical approximation of the gradient $\nabla H$ which 
belongs to the characteristic line bundle of the Hamiltonian flow
\begin{eqnarray}
    \frac{\bZ-\bz}{h} = J\nabla_z H\left( \frac{\bZ+\bz}{2} + hA\left( -J\frac{\bZ-\bz}{h} \right) \right).
    \label{eqn:new}
\end{eqnarray}

In this form, the linear dependency of $\gamma$ and $\beta$ on $h$ has a
geometrical meaning and provides some insight into the development of
an analytical expression for the matrix 
$A$, which gives us the exact flow of the Hamiltonian system. 
This is the subject of a companion article \cite{Jim15c}.

\section{Conclusions and perspectives}
In this paper we have implemented and tested the implicit symplectic
integrators constructed in \cite{Jim15a}. This implementation considers
a simple and coarse iterative process which can be refined.
The main goal here was to test the accuracy of the symplectic
scheme. The construction process developed in \cite{Jim15a}
opens a new point of view on generating functions and symplectic 
integrators owing to several new results. 

$\imath)$ First, we showed that
gauge elements in the primitive 1-form $\theta_{\alpha,\beta,\gamma}$ 
enable the control of the numerical error in the energy integral and other integrals 
by construction. The set of 
parameters $\{\alpha,\beta,\gamma\}=\{\alpha_i,\beta_i,\gamma_i\}$, $i=1,\cdots,n$
with $3n$ elements, gives the set of all possible implicit 
symplectic integrators passing by $\bz$ and $\bZ$, which are consistent with the projection
(\ref{eqn:proj:1}) introduced in \cite{Jim15a}. In the numerical tests, we searched for values
of the parameters that produce the more accurate solution 
for the energy integral and the results are very promising. 

$\imath\imath)$ Second, 
numerical experiments
support the existence of 
exact numerical sympectic integrators, contrary to Ge \cite{Ge91} and
Ge and Marsden's \cite{GM88} claims. This subject is further elaborated in a short note 
\cite{Jim15c} where Ge-Marsden's claim is considered in the case of 
an implicit symplectic integrator. Their claim of the non-existence of 
symplectic integrators exactly preserving the energy integral, relies 
on explicit symplectic maps generating the integrator. 
They also assume an extrapolation by
Taylor series expansion. Our symplectic integrator is implicit, and the 
approximation is done by looking for internal points on the line's flow. 

$\imath\imath\imath)$ Third, 
in the perturbed Hamiltonian pendulum,
the parameters $\beta$ and $\gamma$ and the stepsize $h$ satisfy a linear relationship. This fact,
together with expression
(\ref{eqn:new}) gives a clue on the way the difference equation recovers
information from the Hamiltonian formalism to produce a very accurate numerical solution.
This subject is developed in \cite{Jim15c} with a finer analysis on the 
symplectic transformation which better approximates the flow 
of a Hamiltonian vector field $X_H$.
Other Liouvillian forms as the one obtained from the Poincar\'e's generating 
function are studied in \cite{Jim15f}.

Further work is necessary to understand the subject of 
symplectic integrators using Liouvillian forms. However this series of
papers gives a new point of view on the subject and on the control of 
numerical errors in symplectic integrators.

\section*{Acknowledgements}
This research was developed with support from the Fondation du Coll\`ege de 
France and Total under the research convention PU14150472, as well as the ERC Advanced Grant 
WAVETOMO, RCN 99285, Subpanel PE10 in the F7 framework.

\bibliography{SympInt}
\bibliographystyle{plain}

\end{document}